\documentclass[12pt]{article}
\usepackage{amssymb}
\usepackage{mathrsfs}
\makeatletter
 
 \@addtoreset{equation}{section}
\makeatother

\newcommand{\prf}{\noindent{\bf Proof. }}

\newcommand{\rem}{\noindent{\bf Remark. }}

\newcommand{\ff}[1]{%
    {\bf F}_{#1}}

\newcommand{\qed}{\hbox{\rule[-2pt]{5pt}{11pt}}}

\newtheorem{dfn}{Definition}[section]
\newtheorem{thm}[dfn]{Theorem}
\newtheorem{prop}[dfn]{Proposition}

\newtheorem{lem}[dfn]{Lemma}
\newtheorem{exam}[dfn]{Example}
\newtheorem{conj}[dfn]{Conjecture}

\begin{document}
\title{On some families of divisible formal weight enumerators and their zeta functions}
\author{Koji Chinen\footnotemark[1]}
\date{}
\maketitle

\begin{abstract}
The formal weight enumerators were first introduced by M. Ozeki, and it was shown in the author's previous paper that there are various families of divisible formal weight enumerators. Among them, three families are dealt with in this paper and their properties are investigated: they are analogs of the Mallows-Sloane bound, the extremal property, the Riemann hypothesis, etc. In the course of the investigation, some generalizations of the theory of invariant differential operators developed by I. Duursma and T. Okuda are deduced. 
\end{abstract}

\footnotetext[1]{Department of Mathematics, School of Science and Engineering, Kindai University. 3-4-1, Kowakae, Higashi-Osaka, 577-8502 Japan. E-mail: chinen@math.kindai.ac.jp}
\footnotetext{This work was supported by JSPS KAKENHI Grant Number JP26400028. }


\section{Introduction}\label{section:intro}
The formal weight enumerators were first introduced to coding theory and number theory by Ozeki \cite{Oz}. Recently, the present author \cite{Ch3} showed that there are many other families of ``divisible formal weight enumerators''. So, first we give the definitions of formal weight enumerators and their divisibility. In the following, the action of a matrix $\sigma=\left(\begin{array}{cc} a & b \\ c & d \end{array}\right)$ on a polynomial $f(x,y)\in {\bf C}[x,y]$ is defined by 
\begin{equation}\label{eq:action_sigma}
f^\sigma(x,y)=f(ax+by, cx+dy).
\end{equation}
\begin{dfn}\label{dfn:fwe}
We call a homogeneous polynomial 
\begin{equation}\label{eq:fwe_form}
W(x,y)=x^n+\sum_{i=d}^n A_i x^{n-i} y^i\in {\bf C}[x,y]\quad (A_d\ne 0)
\end{equation}
a formal weight enumerator if 
\begin{equation}\label{eq:fwe_transf}
W^{\sigma_q}(x,y)=-W(x,y)
\end{equation}
for some $q\in {\bf R}$, $q> 0$, $q\ne 1$, where
\begin{equation}\label{eq:macwilliams}
\sigma_q=\frac{1}{\sqrt{q}}\left(\begin{array}{rr} 1 & q-1 \\ 1 & -1 \end{array}\right).
\end{equation}
Moreover, for some fixed $c\in {\bf N}$, we call $W(x,y)$ divisible by $c$ if 
$$A_i\ne 0\quad \Rightarrow \quad c|i$$
is satisfied. 
\end{dfn}
The transformation defined by $\sigma_q$ is often called the MacWilliams transform. Ozeki's formal weight enumerators are of the form
$$W_{{\mathcal H}_{8}}(x,y)^l W_{12}(x,y)^{2m+1}$$
and their suitable linear combinations, where, 
\begin{eqnarray}
W_{{\mathcal H}_{8}}(x,y) &=& x^8+14x^4y^4+y^8, \label{eq:we_hamming}\\
W_{12}(x,y) &=& x^{12}-33x^8y^4-33x^4y^8+y^{12}.\label{eq:w12}
\end{eqnarray}
The polynomial $W_{{\mathcal H}_{8}}(x,y)$ is the weight enumerator of the famous extended Hamming code ${\mathcal H}_{8}$. We have ${W_{{\mathcal H}_{8}}}^{\sigma_2}(x,y)=W_{{\mathcal H}_{8}}(x,y)$ and ${W_{12}}^{\sigma_2}(x,y)=-W_{12}(x,y)$, so Ozeki's formal weight enumerators are those for $q=2$ and $c=4$. 

In the paper \cite{Ch3}, it was shown that the formal weight enumerators divisible by two exist for $q=2, 4, 4/3, 4\pm 2\sqrt{2}, 2\pm 2\sqrt{5}/5, 8\pm 4\sqrt{3}$, etc. The properties of formal weight enumerators vary according to the values of $q$. In this paper, we consider the cases $q=2,4$ and $4/3$. For the cases of other $q$, the reader is referred to \cite{Ch3}. We are mainly interested in the extremal property and the Riemann hypothesis for the zeta functions of  the formal weight enumerators. 

Zeta functions of this kind were first introduced by Duursma \cite{Du1} for the weight enumerators of linear codes, whose theory was developed in his subsequent papers \cite{Du2} -- \cite{Du4}. Later the present author generalized them to Ozeki's formal weight enumerators in \cite{Ch1}, and to some other invariant polynomials in \cite{Ch2}. The definition is the following: 
\begin{dfn}\label{dfn:zeta}
For any homogeneous polynomial of the form (\ref{eq:fwe_form}) and $q\in{\bf R}$ ($q>0, q\ne 1$), there exists a unique polynomial $P(T)\in{\bf C}[T]$ of degree at most $n-d$ such that
\begin{equation}\label{eq:zeta_duursma}
\frac{P(T)}{(1-T)(1-qT)}(y(1-T)+xT)^n=\cdots +\frac{W(x,y)-x^n}{q-1}T^{n-d}+ \cdots.
\end{equation}
We call $P(T)$ and $Z(T)=P(T)/(1-T)(1-qT)$ the zeta polynomial and the zeta function of $W(x,y)$, respectively. 
\end{dfn}
For the proof of existence and uniqueness of $P(T)$, see \cite[Appendix A]{Ch2} for example. Recall that we must assume $d, d^\perp\geq 2$ where $d^\perp$ is defined by 
$$W^{\sigma_q}(x,y)=\pm x^n + A_{d^\perp} x^{n-d^\perp} y^{d^\perp}+ \cdots,$$
when considering the zeta functions (see \cite[p.57]{Du2}). 

If a (formal) weight enumerator $W(x,y)$ has the property $W^{\sigma_q}(x,y)=\pm W(x,y)$, then the zeta polynomial $P(T)$ has the functional equation 
\begin{equation}\label{eq:func_eq}
P(T)=\pm P\left(\frac{1}{qT}\right)q^g T^{2g}\qquad (g=n/2+1-d).
\end{equation}
The quantity $g$ is called the genus of $W(x,y)$. Note that 
\begin{equation}\label{eq:genus_geq0}
d\leq \frac{n}{2}+1
\end{equation}
because $g$ must satisfy $g\geq 0$. Now we can formulate the Riemann hypothesis: 
\begin{dfn}[Riemann hypothesis]\label{dfn:RH}
A (formal) weight enumerator $W(x,y)$ with $W^{\sigma_q}(x,y)=\pm W(x,y)$ satisfies the Riemann hypothesis if all the zeros of $P(T)$ have the same absolute value $1/\sqrt{q}$. 
\end{dfn}
We know examples of (formal) weight enumerators both satisfying and not satisfying the Riemann hypothesis (see \cite[Section 4]{Du3}, \cite{Ch1} -- \cite{Ch4}). 

In the case of the formal weight enumerators treated in this article (especially the cases $q=2$ and $4$), there seems to be similar structures to the cases of the weight enumerators of self-dual codes over the fields $\ff 2$ and $\ff 4$ (so-called Type I and Type IV codes). One of the main purposes of this paper is to investigate such formal weight enumerators and to clarify the properties in common with the weight enumerators of Types I and IV. Our main results are Theorem \ref{thm:analog_mallows-sloane} which establishes analogs of the Mallows-Sloane bound (see Theorem \ref{thm:mallows-sloane}), and Theorem \ref{thm:equiv_rh} which is an analog of Okuda's theorem (see \cite[Theorem 5.1]{Ok}) concerning a certain equivalence of the Riemann hypothesis between some sequences of extremal weight enumerators. 

To this end, we apply the theory of invariant differential operators on invariant polynomial rings, which was introduced by Duursma \cite{Du4} and generalized by Okuda \cite{Ok}. Our second purpose is to generalize their theory further and state it in a form a little easier to use (our main result in this direction is Theorem \ref{thm:duursma_okuda}). 

As to the formal weight enumerators for $q=4/3$, we also find similar structures, but a little different treatment is required. For example, to deduce an analog of the Mallows-Sloane bound (Theorem \ref{thm:mallows-sloane_q=4/3}), it seems that the theory of invariant differential operators does not work well, so we must appeal to the analytical method in MacWilliams-Sloane \cite[p.624-628]{MaSl}. Our main results for this case are Theorem \ref{thm:mallows-sloane_q=4/3} (an analog of the Mallows-Sloane bound) and Theorem \ref{thm:equiv_rh_4over3} (some equivalence of the Riemann hypothesis). 

The rest of the paper is organized as follows: in Section \ref{section:duursma_okuda}, we show the theorem which generalizes the results of Duursma and Okuda. Section \ref{section:fwe_2_4} is devoted to the analysis of divisible formal weight enumerators for $q=2$ and $4$. In Section \ref{section:fwe_4over3}, we discuss the properties of divisible formal weight enumerators for $q=4/3$. 

For a real number $x$, $[x]$ means the greatest integer not exceeding $x$. The Pochhammer symbol $(a)_n$ means $(a)_n=a(a+1)\cdots (a+n-1)$ for $n\geq 1$ and $(a)_0=1$. 
\section{Generalization of the theory of Duursma and Okuda}\label{section:duursma_okuda}
The theory of invariant differential operators on some invariant polynomial rings was introduced by Duurma \cite[Section 2]{Du4}. It considerably simplified the proof of the Mallows-Sloane bound (\cite[Theorem 3]{Du4}). Later a certain generalization is deduced by Okuda \cite[Section 5]{Ok}, which was used to prove a kind of equivalence of the Riemann hypothesis between some sequences of extremal self-dual codes (see \cite[Theorem 5.1 and Section 6]{Ok}). Okuda's idea should be highly appreciated, as well as that of Duursma. 

Their theory must have various applications, in fact one of which is our analysis of formal weight enumerators. In this section, we generalize their theory and give several statements in forms useful for applications. 

We adopt a standard notation as to the action of matrices: for a matrix $\sigma=\left(\begin{array}{cc} a & b \\ c & d \end{array}\right)$ and a pair of variables $(x,y)$, we define 
$$(x,y)^\sigma=(ax+by, cx+dy).$$ 
The action of $\sigma$ on a polynomial $f(x,y)\in {\bf C}[x,y]$ is defined in (\ref{eq:action_sigma}) (these are different from the notation of Duursma \cite{Du4}). For a homogenous polynomial $p(x,y)$, $p(x,y)(D)$ means a differential operator obtained by replacing $x$ by $\partial/\partial x$ and $y$ by $\partial/\partial y$. 
\begin{lem}[Duursma]\label{lem:duursma}
Let $A(x,y)$, $p(x,y)$ be homogenous polynomials in ${\bf C}[x,y]$. Suppose two pairs of variables $(u,v)$ and $(x,y)$ are related by $(u,v)=(x,y)^\sigma$ for a matrix $\sigma=\left(\begin{array}{cc} a & b \\ c & d \end{array}\right)$. Then we have 
$$p^{^t \sigma}(u,v)(D)A(u,v)=p(x,y)(D)A^\sigma(x,y).$$
\end{lem}
\prf This is Duursma \cite[Lemma 1]{Du4}. We state a proof briefly because it is omitted in \cite{Du4}. By the chain rule of differentiation, we have 
$$\frac{\partial}{\partial x}A^\sigma (x,y)=\frac{\partial}{\partial x}A(u,v)
=\frac{\partial A}{\partial u}\frac{\partial u}{\partial x}+\frac{\partial A}{\partial v}\frac{\partial v}{\partial x}.$$
Since $(u,v)=(x,y)^\sigma$, we have $\partial u/\partial x=a$, $\partial v/\partial x=c$. Thus, 
$$\frac{\partial}{\partial x}A(u,v)
=\left( a \frac{\partial }{\partial u} + c \frac{\partial }{\partial v} \right)A(u,v).$$
Similarly we have 
$$\frac{\partial}{\partial y}A(u,v)
=\left( b \frac{\partial }{\partial u} + d \frac{\partial }{\partial v} \right)A(u,v).$$
Therefore we have $(\partial/\partial x, \partial/\partial y)=(\partial/\partial u, \partial/\partial v)^{^t \sigma}$ and generally $p(x,y)(D)=p^{^t \sigma}(u,v)(D)$. \qed

\medskip
\noindent The following proposition is a generalization of the discussion of \cite[pp.108-109]{Du4}: 
\begin{prop}\label{prop:duursma}
Let $a(x,y)$, $A(x,y)$, $p(x,y)$ be homogenous polynomials in ${\bf C}[x,y]$ and suppose $\deg a(x,y) \leq \deg A(x,y)-\deg p(x,y)$. If $a(x,y)|p(x,y)(D)A^\sigma (x,y)$, then we have 
\begin{equation}\label{eq:prop_duursma}
a^{\sigma^{-1}}(x,y)|p^{^t \sigma}(x,y)(D)A(x,y).
\end{equation}
\end{prop}
\prf Let $(u,v)=(x,y)^\sigma$. Then, from Lemma \ref{lem:duursma} and the assumption, we have
$$a(x,y)|p^{^t \sigma}(u,v)(D)A(u,v).$$
Since $a(x,y)=a^{\sigma^{-1}}(u,v)$, we have $a^{\sigma^{-1}}(u,v)|p^{^t \sigma}(u,v)(D)A(u,v)$. This is the same as (\ref{eq:prop_duursma}). \qed

\medskip
\noindent\rem The formula 
$$(x-y)^{d^\perp-1}|((q-1)x-y)(D)A(x,y)$$
which is essentially the same as $(x-y)^{d^\perp-1}|((q-1)x-y)(D)A(x,y)q^{n-k}$ on \cite[p.108]{Du4} is obtained by setting 
$$\sigma=\sigma_q,\quad p(x,y)=y, \quad a(x,y)=y^{d^\perp-1},$$
and the formula
$$(x-\zeta^{-1}y)^{d^\perp-1}|((q-1)x-\zeta y)(D)A(x,y)$$
on \cite[p.109]{Du4} ($x-\zeta y$ on the left hand side seems to be a mistake) is obtained by 
$$\sigma=\left(\begin{array}{cc}1 & 0 \\ 0 & \zeta \end{array}\right),\quad p(x,y)=(q-1)x-y, \quad a(x,y)=(x-y)^{d^\perp-1}.$$
In synthesis of the discussion in Section 2 and Lemma 11 in \cite{Du4}, and Okuda \cite[Proposition 5.4]{Ok}, taking applications to formal weight enumerators into consideration, we obtain the following generalized version of their results: 
\begin{thm}\label{thm:duursma_okuda}
Let $a(x,y)$, $A(x,y)$, $p(x,y)$ be the same as in Proposition \ref{prop:duursma}. We suppose 
\begin{eqnarray*}
p^{^t \sigma}(x,y) &=& c_1 p(x,y),\\
A^\sigma (x,y) &=& c_2 A(x,y)
\end{eqnarray*}
($c_i \in {\bf C}$, $c_i\ne 0$) for a linear transformation $\sigma$. Then we have the following:

\medskip
\noindent(i) 
\begin{equation}\label{eq:duursma_okuda-1}
\{p(x,y)(D)A(x,y)\}^\sigma = \displaystyle\frac{c_2}{c_1} p(x,y)(D)A(x,y).
\end{equation}
\noindent(ii) If $a(x,y)|p(x,y)(D)A(x,y)$, then 
$$a^\sigma(x,y)|p(x,y)(D)A(x,y).$$
Moreover, if $(a(x,y), a^\sigma(x,y))=1$, then
$$a(x,y) a^\sigma(x,y)|p(x,y)(D)A(x,y).$$

\medskip
\noindent(iii) Suppose $a(x,y)|p(x,y)(D)A(x,y)$ and put 
\begin{equation}\label{eq:duursma_okuda-3}
p(x,y)(D)A(x,y)=a(x,y) \tilde a(x,y).
\end{equation}
If $a^\sigma (x,y)=c_3 a(x,y)$ ($c_3 \in {\bf C}$, $c_3\ne 0$), then
\begin{equation}\label{eq:duursma_okuda-3res}
\tilde a^\sigma(x,y)=\frac{c_2}{c_1 c_3}\tilde a(x,y).
\end{equation}
\end{thm}
\prf (i) Let $(u,v)=(x,y)^\sigma$. Then, by Lemma \ref{lem:duursma} and the assumption, we have
$$p(u,v)(D)A(u,v)=\frac{c_2}{c_1}p(x,y)(D)A(x,y).$$
This means (\ref{eq:duursma_okuda-1}). 

\medskip
\noindent(ii) We can prove the former claim by replacing $\sigma$ by $\sigma^{-1}$ in Proposition \ref{prop:duursma} (note that $p^{^t\sigma^{-1}}(x,y)=p(x,y)/c_1$ and $A^{\sigma^{-1}}(x,y)=A(x,y)/c_2$). The latter claim is obvious. 

\medskip
\noindent(iii) Let $\sigma$ act on the both sides of (\ref{eq:duursma_okuda-3}). Then, 
$$\frac{c_2}{c_1}p(x,y)(D)A(x,y)=c_3a(x,y)\tilde a^\sigma (x,y)$$
by (i) and the assumption. Using (\ref{eq:duursma_okuda-3}) again, we get the formula  (\ref{eq:duursma_okuda-3res}). \qed

\medskip
\noindent \rem Okuda \cite[Proposition 5.4]{Ok} is essentially the same as the case where $c_1=c_2=1$ in (i), which was used in the proof of \cite[Theorem 5.1]{Ok}. On the other hand, Duursma \cite[Lemma 11]{Du4} is the case where $c_1=c_2=c_3=1$ for some special $a(x,y)$, $p(x,y)$ and $\sigma$ in (iii). Later we will encounter the cases $c_i=\pm 1$. 
\section{Formal weight enumerators for $q=2$ and $4$}\label{section:fwe_2_4}
In this section, we discuss the properties of formal weight enumerators divisible by two for $q=2$ and $4$. Let
\begin{eqnarray}
\varphi_4(x,y) &=& x^4-6x^2y^2+y^4, \label{eq:fwe_q=2_gen}\\
\varphi_3(x,y) &=& x^3-9xy^2. \label{eq:fwe_q=4_gen}
\end{eqnarray}
Then we can easily see that 
\begin{eqnarray*}
{\varphi_4}^{\sigma_2}(x,y) &=& -\varphi_4(x,y),\\
{\varphi_3}^{\sigma_4}(x,y) &=& -\varphi_3(x,y) 
\end{eqnarray*}
(see also \cite[Section 3]{Ch3}). We can also verify that $W_{2,q}(x,y)=x^2+(q-1)y^2$ satisfies ${W_{2,q}}^{\sigma_q}(x,y)=W_{2,q}(x,y)$ for any $q$. Note that $\varphi_4(x,y)$, $\varphi_3(x,y)$ and $W_{2,q}(x,y)$ are invariant under the action of 
$$\tau=\left(\begin{array}{cc} 1 & 0 \\ 0 & -1 \end{array}\right).$$
We form the following polynomial rings: 
\begin{eqnarray}
R_{\rm I}^- &=& {\bf C}[W_{2,2}(x,y), \varphi_4(x,y)]\label{eq:ring_fwe_q=2},\\
R_{\rm IV}^- &=& {\bf C}[W_{2,4}(x,y), \varphi_3(x,y)]\label{eq:ring_fwe_q=4}.
\end{eqnarray}
These are, so to speak, rings of Type I and Type IV formal weight enumerators, respectively, by analogy with those of Type I and Type IV weight enumerators. Type I weight enumerators are those of self-dual codes over ${\ff 2}$ divisible by two (that is, the weights of all the codewords are divisible by two). The ring of them is
$$R_{\rm I} = {\bf C}[W_{2,2}(x,y), W_{{\mathcal H}_{8}}(x,y)]$$
(see (\ref{eq:we_hamming}) for the definition of $W_{{\mathcal H}_{8}}(x,y)$, see also \cite[p.186]{CoSl} for this ring). Similarly, the Type IV weight enumerators are those of self-dual codes over ${\ff 4}$ divisible by two, whose ring is 
$$R_{\rm IV} = {\bf C}[W_{2,4}(x,y), x^6+45x^2y^4+18y^6]$$
(\cite[p.203]{CoSl}). 

\medskip
\rem The rings $R_{\rm I}^-$ and $R_{\rm IV}^-$ are the invariant polynomial rings of the groups $G_{\rm I}^- = \langle \sigma_2 \tau \sigma_2, \tau \rangle$ and $G_{\rm IV}^- = \langle \sigma_4 \tau \sigma_4, \tau \rangle$, respectively. The group $G_{\rm I}^-$ has order 8 and its Molien series are $\Phi_{\rm I}^- (\lambda)=1/\{(1-\lambda^2)(1-\lambda^4)\}$. The group $G_{\rm IV}^-$ has order 6 and its Molien series are $\Phi_{\rm IV}^- (\lambda)=1/\{(1-\lambda^2)(1-\lambda^3)\}$. 

\medskip
Type I formal weight enumerators are the polynomials $W(x,y)$ of the form (\ref{eq:fwe_form}), given by 
\begin{equation}\label{eq:fwe_typeI}
W_{2,2}(x,y)^l \varphi_4(x,y)^{2m+1}\qquad (l,m \geq 0)
\end{equation}
and their suitable linear combinations (note that we need an odd number of $\varphi_4(x,y)$ to have $W^{\sigma_2}(x,y)=-W(x,y)$). Some examples of such linear combinations will be given in Example \ref{exam:typeI_fwe} later. Similarly, Type IV formal weight enumerators are given by 
\begin{equation}\label{eq:fwe_typeIV}
W_{2,4}(x,y)^l \varphi_3(x,y)^{2m+1}\qquad (l,m \geq 0)
\end{equation}
and their suitable linear combinations (see Example \ref{exam:typeIV_fwe} for an example of such a linear combination). 

Our first goal in this section is Theorem \ref{thm:analog_mallows-sloane}. As a preparation for it, we prove the following proposition, which is an analog of \cite[Lemma 2]{Du4}: 
\begin{prop}\label{prop:divide}
(i) Let $W(x,y)$ be a Type I formal weight enumerator with $d\geq 4$ and let $p(x,y)=xy(x^2-y^2)$. Then we have 
\begin{equation}\label{eq:divide_I}
\{xy(x^2-y^2)\}^{d-3}| p(x,y)(D)W(x,y).
\end{equation}

\medskip
\noindent(ii) Let $W(x,y)$ be a Type IV formal weight enumerator with $d\geq 4$ and let $p(x,y)=y(x^2-9y^2)$. Then we have 
\begin{equation}\label{eq:divide_IV}
\{y(x^2-y^2)\}^{d-3}| p(x,y)(D)W(x,y).
\end{equation}
\end{prop}
\prf (i) It is easy to see that $p^{^t \sigma_2}(x,y)=p(x,y)$ and that $W(x,y)=W(y,x)$ since $W(x,y)$ is invariant under $\sigma_2 \tau \sigma_2$. Moreover, since $W(x,y)$ is of the form (\ref{eq:fwe_form}), we have 
$$p(x,y)(D)W(x,y)=C(x^{n-d-1}y^{d-3}+ \cdots + x^{d-3}y^{n-d-1})$$
for some constant $C$. So we have 
$$(xy)^{d-3}|p(x,y)(D)W(x,y)=-p(x,y)(D)W^{\sigma_2}(x,y)$$
(note that the terms $x^{n-d-1}y^{d-3}$ and $x^{d-3}y^{n-d-1}$ do not disappear when $d\geq 4$ because of the inequality (\ref{eq:genus_geq0})). By Proposition \ref{prop:duursma}, 
$$\{(xy)^{d-3}\}^{{\sigma_2}^{-1}} | -p^{^t \sigma_2}(x,y)(D)W(x,y).$$
Since $\{(xy)^{d-3}\}^{{\sigma_2}^{-1}}=\{(x^2-y^2)/2\}^{d-3}$ and $p^{^t \sigma_2}(x,y)=p(x,y)$, we obtain 
$$(x^2-y^2)^{d-3} | p(x,y)(D)W(x,y).$$
We get (\ref{eq:divide_I}) by Theorem \ref{thm:duursma_okuda} (ii) because $((xy)^{d-3}, (x^2-y^2)^{d-3})=1$. 

\medskip
\noindent(ii) First we note the following: 
$$W^{\sigma_4}(x,y)=-W(x,y), \quad W^\tau(x,y)=W(x,y),$$
$$p^{^t \sigma_4}(x,y)=p(x,y), \quad p^{^t \tau}(x,y)=-p(x,y),$$
$$(y^{d-3})^{{\sigma_4}^{-1}}=\{(x-y)/2\}^{d-3}, \quad \{(x-y)^{d-3}\}^{\tau^{-1}}=(x+y)^{d-3}.$$
Using these, we can prove (\ref{eq:divide_IV}) similarly to (i). \qed

\medskip
\rem As the result of this proposition, we must have $4(d-3)\leq n-4$ for Type I formal weight enumerators with $d\geq 4$, and $3(d-3)\leq n-3$ for Type IV formal weight enumerators with $d\geq 4$. 

\medskip
In the case of Types I and IV weight enumerators, that is the members of $R_{\rm I}$ and $R_{\rm IV}$ of the form (\ref{eq:fwe_form}), the following upper bounds of $d$ by $n$ are known: 
\begin{thm}[Mallows-Sloane]\label{thm:mallows-sloane}\rm 
\begin{eqnarray*}
(\mbox{Type I}) & & d\leq 2\left[ \frac{n}{8} \right] +2, \\
(\mbox{Type IV}) & & d\leq 2\left[ \frac{n}{6} \right] +2.
\end{eqnarray*}
\end{thm}
\prf See \cite[Theorem 3]{Du3} for example. \qed

\medskip
\noindent Our next result is the following: 
\begin{thm}\label{thm:analog_mallows-sloane}
(i) Let $W(x,y)$ be a Type I formal weight enumerator of the form (\ref{eq:fwe_form}). Then we have
$$d\leq 2\left[ \frac{n-4}{8} \right] +2.$$
(ii) Let $W(x,y)$ be a Type IV formal weight enumerator of the form (\ref{eq:fwe_form}). Then we have
$$d\leq 2\left[ \frac{n-3}{6} \right] +2.$$
\end{thm}
\prf (i) We assume $d\geq 4$. Let $p(x,y)=xy(x^2-y^2)$ and $a(x,y)=\{xy(x^2-y^2)\}^{d-3}$. Then we have 
$$p^{^t \sigma_2}(x,y)=p(x,y), \quad p^{^t \tau}(x,y)=-p(x,y),$$
$$W^{\sigma_2}(x,y)=-W(x,y), \quad W^\tau(x,y)=W(x,y),$$
$$a^{\sigma_2}(x,y)=a(x,y), \quad a^\tau (x,y)=-a(x,y)$$
(note that $d$ is even). We apply Theorem \ref{thm:duursma_okuda} (iii). For $\sigma=\sigma_2$, we have $c_1=c_3=1$ and $c_2=-1$, for $\sigma=\tau$, we have $c_1=c_3=-1$, $c_2=1$. So the cofactor $\tilde a(x,y)$ in (\ref{eq:duursma_okuda-3}) satisfies 
$$\tilde a^{\sigma_2}(x,y)=-\tilde a(x,y), \quad \tilde a^{\tau}(x,y)=\tilde a(x,y).$$
Moreover, we can see that $\deg \tilde a(x,y)=n-4d+8$ and $\tilde a(x,y)$ has a term $x^{n-4d+8}$ (see Remark after Proposition \ref{prop:divide}). Hence $\tilde a(x,y)$ is a constant times a Type I formal weight enumerator. Especially, $\tilde a(x,y)$ is divided by $\varphi_4(x,y)=x^4-6x^2y^2+y^4$. This, together with Proposition \ref{prop:divide} (i) yields that 
$$\{xy(x^2-y^2)\}^{d-3}(x^4-6x^2y^2+y^4) | p(x,y)(D)W(x,y).$$
Comparing the degrees on the both sides, we obtain 
$$4(d-3)+4\leq n-4.$$
Putting $d=2d'$ ($d'\in{\bf N}$), we have $d'\leq (n-4)/8+1$. Since $d'$ is an integer, it is equivalent to $d'\leq [(n-4)/8]+1$. The conclusion follows immediately for $d\geq 4$. It also holds for $d=2$. 

\medskip
\noindent(ii) We assume $d\geq 4$. The polynomials $p(x,y)=y(x^2-y^2)$, $W(x,y)$ and $a(x,y)=\{y(x^2-y^2)\}^{d-3}$ satisfy
$$p^{^t \sigma_4}(x,y)=p(x,y), \quad p^{^t \tau}(x,y)=-p(x,y),$$
$$W^{\sigma_4}(x,y)=-W(x,y), \quad W^\tau(x,y)=W(x,y),$$
$$a^{\sigma_4}(x,y)=a(x,y), \quad a^\tau (x,y)=-a(x,y).$$
(note that $d$ is even). We can prove similarly to (i) that 
$$\{y(x^2-y^2)\}^{d-3}(x^3-9xy^2) | p(x,y)(D)W(x,y).$$
We obtain the conclusion by comparing the degrees for $d\geq 4$. It also holds for $d=2$. \qed

\medskip
\rem A similar bound is known for Ozeki's formal weight enumerators which are generated by $W_{{\mathcal H}_{8}}(x,y)$ and $W_{12}(x,y)$ (see (\ref{eq:we_hamming}) and (\ref{eq:w12})), that is, 
$$d\leq 4\left[\frac{n-12}{24}\right]+4$$
(compare this with the Mallows-Sloane bound for Type II weight enumerators $d\leq 4[n/24]+4$, \cite[Theorem 3]{Du3} or \cite[Chapter 19, Theorem 13]{MaSl}). See \cite{Ch1} for details. 

\medskip
Now we can define the notion of extremal formal weight enumerators: 
\begin{dfn}\label{dfn:extremal_fwe}
Let $W(x,y)$ be a Type I or Typr IV formal weight enumerator. We call $W(x,y)$ extremal if the equality holds in Theorem \ref{thm:analog_mallows-sloane}. 
\end{dfn}
We can verify that there exists a unique extremal formal weight enumerator for each degree $n$. 
\begin{exam}\label{exam:typeI_fwe} \rm We collect some examples of Type I formal weight enumerators. 

\medskip
\noindent (1) The extremal formal weight enumerator of degree 12 ($d=4$, note that $a(x,y)=xy(x^2-y^2)$). It coincides with $W_{12}(x,y)$ in (\ref{eq:w12}): 
\begin{eqnarray*}
W_{12}(x,y) &=& \frac{1}{8}\left(9W_{2,2}(x,y)^4 \varphi_4(x,y)-\varphi_4(x,y)^3 \right)\\
&=& x^{12}-33x^8y^4-33x^4y^8+y^{12}.
\end{eqnarray*}
We have
\begin{eqnarray*}
p(x,y)(D)W_{12}(x,y) &=& -6336xy(x^2-y^2)(x^4-6x^2y^2+y^4)\\
  &=& -6336 a(x,y)\varphi_4(x,y).
\end{eqnarray*}
(2) The extremal formal weight enumerator of degree 14 ($d=4$): 
\begin{eqnarray*}
W_{14}(x,y)&:=&\frac{1}{16}(17 W_{2,2}(x,y)^5 \varphi_4(x,y)-W_{2,2}(x,y) \varphi_4(x,y)^3) \\
  &=& x^{14}-26x^{10}y^4-39x^8y^6-39x^6y^8-26x^4y^{10}+y^{14}.
\end{eqnarray*}
We have
$$p(x,y)(D)W_{14}(x,y) = -6240 a(x,y)\varphi_4(x,y)W_{2,2}(x,y).$$
(3) The extremal formal weight enumerator of degree 20 ($d=6$, note that $a(x,y)=\{xy(x^2-y^2)\}^{3}$): 
\begin{eqnarray*}
W_{20}(x,y)&:=&\frac{1}{256}(235 W_{2,2}(x,y)^8 \varphi_4(x,y)+10W_{2,2}(x,y)^4 \varphi_4(x,y)^3
  +11\varphi_4(x,y)^5) \\
  &=& x^{20}-190x^{14}y^6+95x^{12}y^8-836x^{10}y^{10}+95x^8y^{12}-190x^6y^{14}+y^{20}.
\end{eqnarray*}
We have
$$p(x,y)(D)W_{20}(x,y) = -319200 a(x,y)\varphi_4(x,y).$$
(4) An example of a non-extremal formal weight enumerator (degree 20, $d=4$):
\begin{eqnarray*}
W'_{20}(x,y)&:=&\frac{1}{16}(15 W_{2,2}(x,y)^8 \varphi_4(x,y) +\varphi_4(x,y)^5) \\
  &=& x^{20}+5x^{16}y^4-240x^{14}y^6+250x^{12}y^8-1056x^{10}y^{10}\\
  & & +250x^8y^{12}-240x^6y^{14}+5x^4y^{16}+y^{20}.
\end{eqnarray*}
We have
\begin{eqnarray*}
p(x,y)(D)W'_{20}(x,y) &=& 1920 a(x,y)\varphi_4(x,y)\\
 & & \cdot (x^8-238x^6y^2+490x^4y^4-238x^2y^6+y^8).
\end{eqnarray*}
Here the polynomial of degree 8 on the right hand side is equal to
$$\frac{1}{8}(121 \varphi_4(x,y)^2- 113 W_{2,2}(x,y)^4),$$
which is invariant under $\sigma_2$. 
\end{exam}
\begin{exam}\label{exam:typeIV_fwe} \rm We show only one example of the extremal Type IV formal weight enumerator (degree 11, $d=4$): 
\begin{eqnarray*}
W_{11}(x,y) &=& \frac{1}{9}(8 W_{2,4}(x,y)^4 \varphi_3(x,y) + W_{2,4}(x,y)\varphi_3(x,y)^3)\\
  &=& x^{11}-30x^7y^4-336x^5y^6-1035x^3y^8-648xy^{10}.
\end{eqnarray*}
For $p(x,y)=y(y^2-9x^2)$, we have
\begin{eqnarray*}
p(x,y)(D)W_{11}(x,y) &=& -720y(x^2-y^2)(x^3-9xy^2)(x^2+3y^2)\\
  &=& -720 a(x,y)\varphi_3(x,y) W_{2,4}(x,y)
\end{eqnarray*}
where $a(x,y)=\{y(x^2-y^2)\}^{d-3}=y(x^2-y^2)$. 
\end{exam}
Some numerical experiments suggest the following:
\begin{conj}\label{conj:RHtype1_4}
All extremal formal weight enumerators of Types I and IV satisfy the Riemann hypothesis. 
\end{conj}
For the extremal Types I and IV formal weight enumerators, we can also prove analogs of \cite[Theorem 12]{Du4} (the former assertion of it) and \cite[Theorem 19]{Du4}. From Theorem \ref{thm:analog_mallows-sloane}, the degree $n$ can be expressed by $d$ in (\ref{eq:fwe_form}) as follows: 
\begin{eqnarray*}
(\mbox{Type I})& &n=4(d-1)+2v, \quad v=0,1,2,3,\\
(\mbox{Type IV})& &n=3(d-1)+2v, \quad v=0,1,2.
\end{eqnarray*}
Using these parameters, we can prove the following: 
\begin{thm}\label{thm:analog_duursma_th12}
(i) Suppose $d\geq 4$. Then extremal Type I formal weight enumerators $W(x,y)$ satisfy
$$(xy^3-x^3y)(D)W(x,y)=(d-2)_3 (n-d) A_d (x^3y-xy^3)^{d-3}(x^2+y^2)^v (x^4-6x^2y^2+y^4).$$

\medskip
\noindent(ii) Suppose $d\geq 4$. Then extremal Type IV formal weight enumerators $W(x,y)$ satisfy
$$(y^3-9x^2y)(D)W(x,y)=(d-2)_3A_d (x^2y-y^3)^{d-3}(x^2+3y^2)^v (x^3-9xy^2).$$
\end{thm}
\prf We can prove this similarly to \cite[Theorem 12]{Du4}. \qed

\medskip
\noindent We assume $d\geq 4$ and $d$ is even. We put $d-2=m$ ($m\geq 2$, $m$ is even). 
\begin{thm}\label{thm:analog_duursma_th19}
(i) Let $W(x,y)$ be an extremal Type I formal weight enumerator of degree $n=4m+2v+4$ ($m\geq 2$, $m$ is even, $v=0,1,2,3$) and $P(T)=\sum_{i=0}^r p_i T^i$ be the zeta polynomial of $W(x,y)$. Then
\begin{eqnarray}
\sum_{i=0}^{2m+2v+2} p_i {{4m+2v}\choose{m-1+i}}(x-y)^{3m+2v+1-i} y^{m-1+i} 
 &=& \frac{(d-2)_3 (n-d) A_d}{(n-3)_4} (xy)^{m-1} (x^2-y^2)^{m-1}\nonumber\\
 & & \cdot (x^2+y^2)^v (x^4-6x^2y^2+y^4).\label{eq:analog_duursma_th19_I}
\end{eqnarray} 

\medskip
\noindent(ii) Let $W(x,y)$ be an extremal Type IV formal weight enumerator of degree $n=3m+2v+3$ ($m\geq 2$, $m$ is even, $v=0,1,2$) and $Q(T)=P(T)(1+2T)=\sum_{i=0}^r q_i T^i$, where $P(T)$ is the zeta polynomial of $W(x,y)$. Then
\begin{eqnarray}
\sum_{i=0}^{m+2v+2} q_i {{3m+2v}\choose{m-1+i}}(x-y)^{2m+2v+1-i} y^{m-1+i} 
 &=& \frac{(d-2)_3 A_d}{3(n-2)_3} y^{m-1} (x^2-y^2)^{m-1}\nonumber\\
 & & \cdot (x^2+3y^2)^v (x^3-9xy^2).\label{eq:analog_duursma_th19_IV}
\end{eqnarray} 
\end{thm}
\prf Similar to the proof of \cite[Theorem 19]{Du4}. \qed

\medskip
\noindent Unfortunately, we cannot prove Conjecture \ref{conj:RHtype1_4} using Theorem \ref{thm:analog_duursma_th19}. The obstacles are the existence of the factor $x^4-6x^2y^2+y^4$ and $x^3-9xy^2$ on the right hand side of (\ref{eq:analog_duursma_th19_I}) and (\ref{eq:analog_duursma_th19_IV}), as well as $x^{m-1}$ in (\ref{eq:analog_duursma_th19_I}), as was the case of the Type I extremal weight enumerators. However, we can prove a certain equivalence between the Riemann hypothesis for two sequences of extremal formal weight enumerators, which is an analog of Okuda \cite[Theorem 5.1]{Ok}: 
\begin{thm}\label{thm:equiv_rh}
(i) Let $W(x,y)$ be the extremal Type I formal weight enumerator of degree $n=8k+4$ ($k\geq 1$) with the zeta polynomial $P(T)$. Then 
$$W^\ast(x,y):=\frac{1}{n(n-1)}(x^2+y^2)(D)W(x,y)$$
is the extremal formal weight enumerator of degree $8k+2$ with the zeta polynomial $(2T^2-2T+1)P(T)$. The Riemann hypothesis for $W(x,y)$ is equivalent to that of $W^\ast(x,y)$. 

\medskip
\noindent (ii) Let $W(x,y)$ be the extremal Type IV formal weight enumerator of degree $n=6k+3$ ($k\geq 1$) with the zeta polynomial $P(T)$. Then 
$$W^\ast(x,y):=\frac{1}{n(n-1)}\left(x^2+\frac{1}{3}y^2\right)(D)W(x,y)$$
is the extremal formal weight enumerator of degree $6k+1$ with the zeta polynomial $(4T^2-2T+1)P(T)/3$. The Riemann hypothesis for $W(x,y)$ is equivalent to that of $W^\ast(x,y)$. 
\end{thm}
\prf (i) We follow the method of Okuda \cite[Section 5]{Ok}. Our proof is similar to it, but we state a proof because \cite{Ok}, being written in Japanese, is not easily accessible to all the readers. We have $W^{\sigma_2}(x,y)=-W(x,y)$ and $W^\tau(x,y)=W(x,y)$. For $p(x,y)=x^2+y^2$, we have $p^{^t\sigma_2}(x,y)=p^{^t\tau}(x,y)=p(x,y)$. So, from Theorem \ref{thm:duursma_okuda} (i) (the case $c_1=1$, $c_2=-1$), we can see that $W^\ast(x,y)$ is a formal weight enumerator of degree $n-2$, the term of smallest degree with respect to $y$ is that of $x^{n-d} y^{d-2}$. If $n=8k+4$, then $2[(n-4)/8]+2=2k+2$, and if $n=8k+2$, then $2[(n-4)/8]+2=2k$. Since the extremal formal weight enumerator is determined uniquely for each degree $n$, we can see that $W^\ast(x,y)$ is extremal. 

To deduce the relation between the zeta polynomials, we need the MDS weight enumerators for $q=2$. Let $M_{n,d}=M_{n,d}(x,y)$ be the $[n, k=n-d+1, d]$ MDS weight enumerator and suppose the genus of $W(x,y)$ is $n/2+1-d$. Then $P(T)=\sum_{i=0}^{n-2d+2} a_i T^i$ is related to $W(x,y)$ by 
\begin{equation}\label{eq:rel_P(T)_W(x,y)}
W(x,y)=a_0 M_{n,d} + a_1 M_{n,d+1} + \cdots + a_{n-2d+2}M_{n,n-d+2}
\end{equation}
(see \cite[formula (5)]{Du2}). Note that $d\geq 4$. We have
\begin{eqnarray*}
x(D)M_{n,i}(x,y) &=& n M_{n-1,i}(x,y),\\
y(D)M_{n,i}(x,y) &=& n (M_{n-1,i-1}(x,y)-M_{n-1,i}(x,y))
\end{eqnarray*}
(see ``puncturing and averaging operator'' and ``shortening and averaging operator'' of \cite[Section 3]{Du2}). We act $x(D)$ on both sides of (\ref{eq:rel_P(T)_W(x,y)}) and obtain 
$$x(D)W(x,y)=n(a_0 M_{n-1,d} + a_1 M_{n-1,d+1} + \cdots + a_{n-2d+2} M_{n-1,n-d+2}).$$
So we see that the zeta polynomial of $x(D)W(x,y)/n$ is $P(T)$. Acting $x(D)$ once again, we can see the zeta polynomial of $x^2(D)W(x,y)/n(n-1)$ is $P(T)$, too. For the operator $y(D)$, we have
\begin{eqnarray*}
\frac{1}{n}y(D)W(x,y) &=& a_0 M_{n-1,d-1} + (a_1-a_0) M_{n-1,d}+\cdots + 
   (a_{n-2d+2}-a_{n-2d+1})M_{n-1,n-d+1}\\
   & & -a_{n-2d+2}M_{n-1,n-d+2},
\end{eqnarray*}
of which the zeta polynomial is 
$$a_0+(a_1-a_0)T+\cdots +(a_{n-2d+2}-a_{n-2d+1})T^{n-2d+2}-a_{n-2d+2}T^{n-2d+3}=(1-T)P(T).$$
From this, we can also see that the zeta polynomial of $y^2(D) W(x,y)/n(n-1)$ is $(1-T)^2P(T)$. Note that $x^2(D)W(x,y)/n(n-1)$ begins with the term of $M_{n-2,d}$, whereas $y^2(D) W(x,y)/n(n-1)$ begins with $M_{n-2,d-2}$. Therefore, adjusting the degree, we can conclude that the zeta polynomial of $W^\ast(x,y)$ is 
$$T^2 P(T)+(1-T)^2 P(T)=(2T^2-2T+1)P(T).$$
The equivalence of the Riemann hypothesis is immediate since both roots of $2T^2-2T+1$ have the same absolute value $1/\sqrt{2}$. 

\medskip
\noindent(ii) We use $p(x,y)=x^2+y^2/3$. The proof is similar to that of (i) (this case is almost the same as \cite[Theorem 5.1]{Ok}). \qed
\section{Formal weight enumerators for $q=4/3$}\label{section:fwe_4over3}
In our previous paper \cite{Ch3}, we have found that 
\begin{equation}\label{eq:gen_fwe4over3_deg6}
\varphi_6(x,y) = x^6-5x^4y^2+\frac{5}{3}x^2y^4-\frac{1}{27}y^6
\end{equation}
satisfies ${\varphi_6}^{\sigma_{4/3}}(x,y)=-\varphi_6(x,y)$. We also know that 
\begin{equation}\label{eq:gen_fwe4over3_deg2}
W_{2,4/3}(x,y) = x^2+\frac{1}{3}y^2
\end{equation}
satisfies ${W_{2,4/3}}^{\sigma_{4/3}}(x,y)=W_{2,4/3}(x,y)$. So we form the following two polynomial rings
\begin{eqnarray*}
R_{4/3}^- &=&{\bf C}[\varphi_6(x,y), W_{2,4/3}(x,y)],\\
R_{4/3} &=&{\bf C}[\varphi_6(x,y)^2, W_{2,4/3}(x,y)]. 
\end{eqnarray*}
The formal weight enumerators are polynomials of the form (\ref{eq:fwe_form}) in $R_{4/3}^-$ given by 
$$W_{2,4/3}(x,y)^l \varphi_6(x,y)^{2m+1} \qquad (l,m\geq 0)$$
and their suitable linear combinations. We also consider the invariant polynomials of the form (\ref{eq:fwe_form}) in $R_{4/3}$ given by 
$$W_{2,4/3}(x,y)^l \varphi_6(x,y)^{2m} \qquad (l,m\geq 0, \quad (l,m)\ne(0,0))$$
and their suitable linear combinations. 

We show that the rings $R_{4/3}^-$ and $R_{4/3}$ can be realized as invariant polynomial rings of some groups in $SL_2({\bf C})$. We can see that $R_{4/3}$ is indeed the largest ring which contains polynomials invariant under $\sigma_{4/3}$ and divisible by two. We showed in \cite{Ch3} that there is no $W(x,y)$ of degree less than six satisfying $W^{\sigma_{4/3}}(x,y)=-W(x,y)$, so $R_{4/3}^-$ is also the largest ring of formal weight enumerators for $q=4/3$ divisible by two. 
\begin{prop}\label{prop:groups_4over3}
(i) Let $\eta=\displaystyle\frac{1}{2}\left(\begin{array}{rr} 1 & 1 \\ -3 & 1 \end{array}\right)$, 
$\tau=\left(\begin{array}{rr} 1 & 0 \\ 0 & -1 \end{array}\right)$ and $G_{4/3}^-=\langle \eta, \tau \rangle$. Then we have $|G_{4/3}^-|=12$ and the Molien series are
$$\Phi(\lambda)=\frac{1}{(1-\lambda^2)(1-\lambda^6)}.$$
The ring $R_{4/3}^-$ is the invariant polynomial ring of $G_{4/3}^-$. 

\medskip
\noindent(ii) Let $G_{4/3}=\langle \sigma_{4/3}, \tau \rangle$. Then we have $|G_{4/3}|=24$ and the Molien series are
$$\Phi(\lambda)=\frac{1}{(1-\lambda^2)(1-\lambda^{12})}.$$
The ring $R_{4/3}$ is the invariant polynomial ring of $G_{4/3}$. 
\end{prop}
\prf (i) We can verify that $\eta$ has order 6, $\tau^2=I$ ($I$ is the identity matrix) and the relation $\tau\eta=\eta^5\tau$. It follows that
\begin{eqnarray*}
G_{4/3}^- &=& \{\eta^i \tau^j\ ;\ 0\leq i \leq 5, \ j=0,1 \} \\
    &=&\langle \eta \rangle \rtimes \langle \tau \rangle. 
\end{eqnarray*}
Thus we can see $|G_{4/3}^-|=12$. The Molien series can be calculated directly by the definition (\cite[p.600]{MaSl})
$$\Phi(\lambda)=\frac{1}{|G_{4/3}^-|}\sum_{A\in G_{4/3}^-} \frac{1}{\det(I-\lambda A)}.$$
The result implies that the invariant polynomial ring ${\bf C}[x,y]^{G_{4/3}^-}$ has two generators, one of which has degree two and the other has degree six. It can be checked that $\eta$ and $\tau$ fix both $W_{2,4/3}(x,y)$ and $\varphi_6(x,y)$. 

\medskip
\noindent(ii) We have ${\sigma_{4/3}}^2=\rho^2=I$, $\sigma_{4/3}\tau$ has order 12 and so $\tau\sigma_{4/3}=(\sigma_{4/3}\tau)^{11}$. There are no $k,l\in{\bf Z}$ such that $(\sigma_{4/3}\tau)^k\sigma_{4/3}=(\sigma_{4/3}\tau)^l$. Therefore
\begin{eqnarray*}
G_{4/3} &=& \{(\sigma_{4/3}\tau)^i {\sigma_{4/3}}^j\ ;\ 0\leq i \leq 11, \ j=0,1 \} \\
    &=&\langle \sigma_{4/3}\tau \rangle \rtimes \langle \sigma_{4/3} \rangle 
\end{eqnarray*}
and $|G_{4/3}|=24$. The Molien series are obtained similarly. It is obvious that ${\sigma_{4/3}}$ and $\tau$ fix $\varphi_6(x,y)^2$ and $W_{2,4/3}(x,y)$. \qed

\medskip
Next we consider analogs of Mallows-Slaone bound. In the present case ($q=4/3$), it seems difficult to find a good differential operator $p(x,y)(D)$ and a good polynomial $a(x,y)$ like in the previous section, but it is possible to prove the following by use of an analytic method of \cite[Chapter 19, Section 5]{MaSl}: 
\begin{thm}\label{thm:mallows-sloane_q=4/3}
(i) Any invariant polynomial of the form (\ref{eq:fwe_form}) in $R_{4/3}$ satisfies 
\begin{equation}\label{eq:mallows-sloane_q=4/3-1}
d\leq 2\left[\frac{n}{12}\right]+2.
\end{equation}

\medskip
\noindent(ii) Any formal weight enumerator of the form (\ref{eq:fwe_form}) in $R_{4/3}^-$ with $n\equiv 6\ ({\rm mod}\ 12)$ satisfies 
\begin{equation}\label{eq:mallows-sloane_q=4/3-2}
d\leq 2\left[\frac{n-6}{12}\right]+2.
\end{equation}
\end{thm}
\prf (i) We follow the method of \cite[p.624-628]{MaSl}. So we use a similar notation and state an outline. Let 
$$W_2(x,y)=W_{2,4/3}(x,y)$$
and 
\begin{eqnarray*}
W'_{12}(x,y) &=& \frac{1}{2}(W_{2,4/3}(x,y)^6-\varphi_6(x,y)^2)\\
   &=&\frac{1}{81}x^2y^2(x^2-y^2)^2(9x^2-y^2)^2.
\end{eqnarray*}
Then we have
$$R_{4/3}={\bf C}[W_2(x,y), W'_{12}(x,y)].$$
An invariant polynomial $W(x,y)$ in $R_{4/3}$ of the form (\ref{eq:fwe_form}) can be written as
\begin{equation}\label{eq:lin_comb}
W(x,y)=\sum_{r=0}^\mu a_r W_2(x,y)^{6\mu+\nu-6r} W'_{12}(x,y)^r,
\end{equation}
here, $n=\deg W(x,y)=2(6\mu+\nu)$ ($\mu\geq 0$, $0\leq \nu \leq 5$, $(\mu,\nu)\ne(0,0)$). Suppose we choose suitable $a_r$ and we cancel as many coefficients as possible. The right hand side of (\ref{eq:lin_comb}) is a linear combination of $\mu + 1$ polynomials, so we can at least make $y^2, y^4, \cdots, y^{2\mu}$ disappear. So we assume
\begin{equation}\label{eq:lin_comb_cancel}
W(x,y)=x^n + \sum_{r=\mu+1}^{6\mu+\nu} A_{2r} x^{n-2r} y^{2r}.
\end{equation}
Our goal is to prove $A_{2\mu+2}\ne0$. We substitute $x$ by 1 and $y^2$ by $x$ in $W_2(x,y)$ and $W'_{12}(x,y)$. We put
\begin{eqnarray*}
f(x) &=& 1+\frac{1}{3}x,\\
g(x) &=& x(1-x)^2 (1-x/9)^2.
\end{eqnarray*}
The function ${\mit\Phi}(x)=xf(x)^6/g(x)$ satisfies the conditions of the B\"urmann-Lagrange Theorem (see \cite[Chapter 19, Theorem 14]{MaSl}) and we can conclude that 
\begin{equation}\label{eq:coeff_A_2m+2}
A_{2\mu+2} = \frac{9^{2\mu+2}(6\mu+\nu)}{3\cdot(\mu+1)!} \frac{d^\mu}{dx^\mu} 
              \left.\left\{ \frac{(1+x/3)^{5-\nu}}{(x-1)^{2\mu+2}(x-9)^{2\mu+2}} \right\}\right|_{x=0}.
\end{equation}
Let 
$$F_\mu (x;\alpha, \beta)=(x-\alpha)^{-2\mu-2}(x-\beta)^{-2\mu-2}$$
for $\alpha, \beta>0$. Then it is easy to see that 
$$F_\mu^{(l)}(0;\alpha, \beta)=\sum_{r=0}^l {{l}\choose{r}} (2\mu+2)_{l-r}(2\mu+2)_r \alpha^{-2\mu-2-l+r} \beta^{-2\mu-2-r}>0$$
for all $l\geq 0$ ($\alpha=1$, $\beta=9$ in our case). Moreover, since $5-\nu\geq 0$, we have  $\{(1+x/3)^{5-\nu}\}^{(l)}|_{x=0}>0$ unless $\{(1+x/3)^{5-\nu}\}^{(l)}$ is identically zero. Thus we can see that $A_{2\mu+2}>0$ for all $\mu\geq 0$ and that $d\leq 2\mu+2$. We recall $n=2(6\mu+\nu)$ and $d$ is even. Putting $d=2d'$ ($d'\in{\bf N}$), we obtain
$$d'\leq \mu+1 =\frac{n}{12}-\frac{\nu}{6}+1\leq \frac{n}{12}+1, \quad d'\leq \left[\frac{n}{12}\right]+1.$$
The conclusion follows immediately. 

\medskip
\noindent(ii) The proof is similar to (i), but a little more delicate estimate is needed. If we cancel as many coefficients as possible, the formal weight enumerator $W(x,y)$ of degree $n\equiv 6\ ({\rm mod}\ 12)$ can be written in the form 
\begin{eqnarray*}
W(x,y) &=& \sum_{r=0}^\mu b_r W'_{12}(x,y)^r \varphi_6(x,y)^{2\mu-2r+1} \qquad (\mu\geq 0)\nonumber\\
  &=& x^{12\mu+6} + \sum_{r=\mu+1}^{6\mu+3} A_{2r} x^{12\mu-2r+6} y^{2r}. 
\end{eqnarray*}
Here, $n=\deg W(x,y)=12\mu+6$. We put 
\begin{eqnarray*}
f(x) &=& 1-5x+\frac{5}{3}x^2-\frac{1}{27}x^3,\\
g(x) &=& x(1-x)^2 (1-x/9)^2.
\end{eqnarray*}
By a similar argument to (i), we get 
\begin{eqnarray}
A_{2\mu+2} &=& -\frac{9^{2\mu+2}(2\mu+1)}{(\mu+1)!}\frac{d^\mu}{dx^\mu} 
\left.\left\{\left(-\frac{1}{9}x^2+\frac{10}{3}x-5\right)F_\mu(x; 1, 9)\right\}\right|_{x=0}\nonumber\\
 &=& -\frac{9^{2\mu+2}(2\mu+1)}{(\mu+1)!}
     \left\{5\sum_{r=0}^\mu {{\mu}\choose{r}} (2\mu+2)_{\mu-r}(2\mu+2)_{r}9^{-2\mu-2-r} \right. \nonumber\\
 & &-\frac{10}{3}\mu\sum_{r=0}^{\mu-1} {{\mu-1}\choose{r}} (2\mu+2)_{\mu-1-r}(2\mu+2)_{r}9^{-2\mu-2-r} \nonumber\\
 & &+\left.\frac{\mu(\mu-1)}{9}\sum_{r=0}^{\mu-2} {{\mu-2}\choose{r}} (2\mu+2)_{\mu-2-r}(2\mu+2)_{r}9^{-2\mu-2-r}
\right\}\label{eq:coeff_A_2m+2_fwe}
\end{eqnarray}
for $\mu\geq2$. Now we prove $A_{2\mu+2}<0$. It suffices to show that 
$$\sum_{r=0}^\mu {{\mu}\choose{r}} (2\mu+2)_{\mu-r}(2\mu+2)_{r}9^{-2\mu-2-r}
>\mu\sum_{r=0}^{\mu-1} {{\mu-1}\choose{r}} (2\mu+2)_{\mu-1-r}(2\mu+2)_{r}9^{-2\mu-2-r},$$
that is, to show that 
$$\sum_{r=0}^{\mu-1} {{\mu}\choose{r}} (2\mu+2)_{\mu-r}(2\mu+2)_{r}9^{-2\mu-2-r}+
{{\mu}\choose{\mu}} (2\mu+2)_\mu 9^{-3\mu-2}$$
\begin{equation}\label{eq:ineq_goal}
>\mu\sum_{r=0}^{\mu-1} {{\mu-1}\choose{r}} (2\mu+2)_{\mu-1-r}(2\mu+2)_{r}9^{-2\mu-2-r}.
\end{equation}
If $r\ne0$, then ${{\mu}\choose{r}}> {{\mu-1}\choose{r}}$. If $0\leq r \leq \mu-1$, then we have
$$\frac{(2\mu+2)_{\mu-r}(2\mu+2)_{r}}{\mu (2\mu+2)_{\mu-1-r}(2\mu+2)_{r}}
=\frac{3\mu+1-r}{\mu}>\frac{2\mu+2}{\mu}>1.$$
From these, we can prove (\ref{eq:ineq_goal}) and get $A_{2\mu+2}<0$. Since $n=12\mu+6$, we can estimate $d$ as 
$$d\leq 2\mu+2=2\cdot\frac{n-6}{12}+2,$$
the conclusion follows similarly to (i) for $\mu\geq 2$, that is, $n\geq 30$. For the cases $\mu=0,1$, explicit constructions show the bound: when $\mu=0$ ($n=6$), there is only one formal weight enumerator $\varphi_6(x,y)$ whose $d=2$, so (\ref{eq:mallows-sloane_q=4/3-2}) holds. When $\mu=1$ ($n=18$), the basis contains two formal weight enumerators $\varphi_6(x,y)^3$ and $W'_{12}(x,y)\varphi_6(x,y)$. We eliminate the term of $y^2$ by making
\begin{eqnarray*}
\varphi_6(x,y)^3 + 15 W'_{12}(x,y)\varphi_6(x,y) 
&=& x^{18} - \frac{85}{3}x^{14}y^4 + \frac{1037}{27}x^{12}y^6 - \frac{935}{27}x^{10}y^8 \\
 & & + \frac{935}{81}x^8y^{10} - \frac{1037}{729}x^6y^{12} + \frac{85}{729}x^4y^{18} - \frac{1}{19683}y^{18}
\end{eqnarray*}
whose $d=4$. Thus we have proved the theorem. \qed
\begin{exam}\label{exam:coeff_A_d}\rm (i) Let $\mu=1$ and $\nu=5$ in (\ref{eq:coeff_A_2m+2}). Then (\ref{eq:coeff_A_2m+2}) gives $A_4$ for $n=\deg W(x,y)=22$: 
$$A_4=\frac{9^4\cdot 11}{3\cdot 2}\frac{d}{dx}\left.\left\{ \frac{1}{(x-1)^4(x-9)^4} \right\}\right|_{x=0}
=\frac{220}{27}.$$
It coincides with the relevant coefficient in 

\medskip
$\displaystyle\frac{1}{36}\{ 25 W_{2,4/3}(x,y)^{11} +11 W_{2,4/3}(x,y)^5 \varphi_6(x,y)^2 \}$
\begin{eqnarray*}
&=&x^{22}+\frac{220}{27}x^{18}y^4+\frac{2497}{243}x^{16}y^6+\frac{2750}{729}x^{14}y^8+\frac{484}{2187}x^{12}y^{10}+\frac{484}{6561}x^{10}y^{12}\\
& &+\frac{2750}{19683}x^8y^{14}+\frac{2497}{59049}x^6y^{16}+\frac{220}{59049}x^4y^{18}+\frac{1}{177147}y^{22}.
\end{eqnarray*}

\medskip
\noindent (ii) Let $\mu=2$ in (\ref{eq:coeff_A_2m+2_fwe}). Then (\ref{eq:coeff_A_2m+2_fwe}) gives 
$$A_6=-\frac{14065}{81}$$
for $n=\deg W(x,y)=30$. It coincides with the relevant coefficient in 

\medskip
$\displaystyle\frac{1}{8424}\{ 10075 W_{2,4/3}(x,y)^{12} \varphi_6(x,y) - 2600 W_{2,4/3}(x,y)^6 \varphi_6(x,y)^3 
+949 \varphi_6(x,y)^5 \}$
$$= x^{30}-\frac{14065}{81}x^{24}y^6+ \cdots .$$
\end{exam}
\medskip
\rem (i) It is very plausible that the bound (\ref{eq:mallows-sloane_q=4/3-2}) holds for any formal weight enumerators in $R_{4/3}^-$. The general case requires the analysis of 
$$\sum_{r=0}^\mu b_r W_2(x,y)^c W'_{12}(x,y)^r \varphi_6(x,y)^{2\mu-2r+1} \quad(0\leq c \leq 5, \mu\geq 0),$$
which is attended with much difficulty. What is treated in Theorem \ref{thm:mallows-sloane_q=4/3} (ii) is the case where $c=0$. 

\medskip
\noindent (ii) One is tempted to find suitable $p(x,y)$ to prove Theorem \ref{thm:mallows-sloane_q=4/3} like in the previous section. One of the candidates of $p(x,y)$ should be
$$p(x,y)=xy(x^2-y^2)(x^2-9y^2)$$
which satisfies $p^{^t\sigma_{4/3}}(x,y)=p(x,y)$ and $p^{^t\tau}(x,y)=-p(x,y)$. Using this and a similar reasoning to the previous section, we can prove 
$$\{xy(x^2-y^2)(9x^2-y^2)\}^{d-5} \varphi_6(x,y) | p(x,y)(D)W(x,y)$$
for a formal weight enumerator $W(x,y)$ in $R_{4/3}^-$ with $d\geq 6$, but this does not reach the desired bound (\ref{eq:mallows-sloane_q=4/3-2}). 

\medskip
\noindent We can define the extremal polynomials in $R_{4/3}$: 
\begin{dfn}\label{dfn:extremal_4over3}
Let $W(x,y)$ be a polynomial of the form (\ref{eq:fwe_form}) in $R_{4/3}$. We call $W(x,y)$ extremal if the equality holds in (\ref{eq:mallows-sloane_q=4/3-1}). 
\end{dfn}
Some numerical experiments suggest the following:
\begin{conj}\label{conj:RH4over3}
All extremal polynomial of the form (\ref{eq:fwe_form}) in $R_{4/3}$ satisfy the Riemann hypothesis. 
\end{conj}
We cannot prove the above conjecture, but we can prove the following theorem, analogous to Theorem \ref{thm:equiv_rh}: 
\begin{thm}\label{thm:equiv_rh_4over3}
Let $W(x,y)$ be the extremal polynomial of the form (\ref{eq:fwe_form}) in $R_{4/3}$ and of degree $n=12k$ ($k\geq 1$) with the zeta polynomial $P(T)$. Then 
$$W^\ast(x,y):=\frac{1}{n(n-1)}(x^2+3y^2)(D)W(x,y)$$
is the extremal polynomial of degree $12k-2$ with the zeta polynomial $(4T^2-6T+3)P(T)$. The Riemann hypothesis for $W(x,y)$ is equivalent to that of $W^\ast(x,y)$. 
\end{thm}
\prf We use $p(x,y)=x^2+3y^2$. We can prove the theorem similarly to Theorem \ref{thm:equiv_rh} (we omit the detail). \qed 
\begin{exam}\label{exam:equiv_rh_4over3}\rm 
The case $k=1$. The extremal polynomial of degree 12 is 
\begin{eqnarray*}
W_{12}^{\rm E}(x,y) &=& \frac{1}{6}\{5 W_{2,4/3}(x,y)^6+\varphi_6(x,y)^2\}\\
  &=& x^{12}+\frac{55}{9}x^8y^4 - \frac{176}{81}x^6y^6 + \frac{55}{81}x^4y^8 + \frac{1}{729}y^{12}.
\end{eqnarray*}
The zeta polynomial is 
$$P_{12}^{\rm E}(T)=\frac{1}{5103}(448T^6+896T^5+1128T^4+1092T^3+846T^2+504T+189).$$
On the other hand, 
\begin{eqnarray*}
(W_{12}^{\rm E})^\ast(x,y) &=& x^{10}+\frac{5}{3}x^8y^2+\frac{10}{9}x^6y^4+\frac{10}{27}x^4y^6+\frac{5}{81}x^2y^8+\frac{1}{243}y^{10}\\
  &=& W_{2,4/3}(x,y)^5,
\end{eqnarray*}
which is indeed the extremal polynomial of degree 10. We can verify that its zeta polynomial coincides with $(4T^2-6T+3)P_{12}^{\rm E}(T)$. 
\end{exam}
\rem It can be conjectured that a theorem similar to Theorem \ref{thm:equiv_rh_4over3} holds for extremal formal weight enumerators in $R_{4/3}^-$. In this case, the relevant degrees are $n=12k+6$ and $12k+4$ ($k\geq 1$). We proved (\ref{eq:mallows-sloane_q=4/3-2}) for the degree $n=12k+6$, but not for the degree $12k+4$. The author observed that there was a relation $P_{28}^{\rm E}(T)=(4T^2-6T+3)P_{30}^{\rm E}(T)$, where $P_{30}^{\rm E}(T)$ is the zeta polynomial of the extremal formal weight enumerator of degree 30 ($d=6$) and $P_{28}^{\rm E}(T)$ is that of the unique formal weight enumerator of degree 28 with $d=4$ (at this degree, we can verify that it is extremal). 

\bigskip
\noindent{\it Acknowledgement. }The author would like to express his sincere gratitude to Professor Yoshio Mimura for his pieces of valuable advice on the structures of matrix groups. This work was established mainly during the author's stay at University of Strasbourg for the overseas research program of Kindai University. He would also like to thank Professor Yann Bugeaud at University of Strasbourg for his hospitality and Kindai University for giving him a chance of the program.


\begin{thebibliography}{99}
\bibitem[1] {Ch1} K. Chinen, Zeta functions for formal weight enumerators and the extremal property, Proc. Japan Acad. 81 Ser. A. (2005), 168-173. 
\bibitem[2] {Ch2} K. Chinen, An abundance of invariant polynomials satisfying the Riemann hypothesis, Discrete Math. 308 (2008), 6426-6440. 
\bibitem[3] {Ch3} K. Chinen, Construction of divisible formal weight enumerators and extremal polynomials not satisfying the Riemann hypothesis, preprint. 
\bibitem[4] {Ch4} K. Chinen, Extremal invariant polynomials not satisfying the Riemann hypothesis, preprint. 

\bibitem[5] {CoSl} J. H. Conway, N. J. A. Sloane, Sphere Packings, Lattices and Groups, third ed., Springer Verlag, NewYork, 1999.

\bibitem[6] {Du1} I. Duursma, Weight distribution of geometric Goppa codes, Trans. Amer. Math. Soc. 351, No.9 (1999), 3609-3639.
\bibitem[7] {Du2} I. Duursma, From weight enumerators to zeta functions, Discrete Appl. Math. 111 (2001), 55-73.
\bibitem[8] {Du3} I. Duursma, A Riemann hypothesis analogue for self-dual codes, DIMACS series in Discrete Math. and Theoretical Computer Science 56 (2001), 115-124.
\bibitem[9] {Du4}I. Duursma, Extremal weight enumerators and ultraspherical polynomials, Discrete Math. 268, No.1-3 (2003), 103-127. 
\bibitem[10] {MaSl} F. J. MacWilliams, N. J. A. Sloane, The Theory of Error-Correcting Codes, North-Holland, Amsterdam, 1977.  
\bibitem[11] {Ok} T. Okuda, Relation between zeta polynomials and differential operators on some invariant rings (in Japanese), RIMS K\^oky\^uroku Bessatsu B20 (2010), 57-69. See also RIMS K\^oky\^uroku 1593 (2008), 145-153. 
\bibitem[12] {Oz} M. Ozeki, On the notion of Jacobi polynomials for codes, Math. Proc. Camb. Phil. Soc. 121 (1997), 15-30. 
\bibitem[13] {ShTo} Shephard, G. C. and Todd, J. A. : Finite unitary reflection groups, Canad. J. Math. {\bf 6} (1954), 274-304. 
\end{thebibliography}
\end{document}